\documentclass[11pt]{article}
\newcommand{\eqdef}{\, =\kern -15pt\raise 6pt\hbox{{\tiny\textrm{def}}}\,}
\newtheorem{theorem}{Theorem}
\begin{document}
\begin{center}{\Huge Patterns and Fractions}

\vspace{1in}

{\Large Aaron Robertson}\\
{\large Temple University, Philadelphia, PA 19122}\\
\texttt{<aaron@euclid.math.temple.edu>}

\vspace{.2in}

{\Large Herbert S. Wilf}\\{\large University of Pennsylvania, Philadelphia,
PA 19104-6395}\\\texttt{<wilf@math.upenn.edu>}

\vspace{.2in}
and
\vspace{.2in}

{\Large Doron Zeilberger}\\{\large Temple University, Philadelphia, PA
19122}\\\texttt{<zeilberg@euclid.math.temple.edu>}

\vspace{.7in}

April 11, 1999
\end{center}
\vspace{1in}
\begin{abstract}
We find, in the form of a continued fraction, the generating function for
the number of $(132)$-avoiding permutations that have a given number of
$(123)$ patterns, and show how to extend this to permutations that have
exactly one $(132)$ pattern. We find some properties of the continued
fraction, which is similar to, though more general than, those that were
studied by Ramanujan.\end{abstract}

\newpage

A $(132)$ pattern (\textit{resp.} a $(123)$ pattern) in a permutation $\pi$
of $|\pi|$ letters is a triple $1\le i<j<k\le n$ of indices for which
$\pi(i)<\pi(k)<\pi(j)$ (\textit{resp.} $\pi(i)<\pi(j)<\pi(k)$). Let
$f_r(n)$ denote the number of permutations $\pi$ of $n$ letters that have
no $(132)$ patterns and exactly $r$ $(123)$ patterns. Our main result is
the following.
\begin{theorem} The generating function for the $\{f_r(n)\}$ is
\begin{equation}
\label{eq:argen}
\sum_{r,n\ge 0}f_r(n)z^nq^r
={1\over\displaystyle 1-{\strut z\over\displaystyle 1-{\strut
z\over\displaystyle 1-{\strut zq\over\displaystyle 1-{\strut
zq^3\over\displaystyle 1-{\strut zq^6\over\displaystyle
\dots}}}}}}\end{equation}
in which the $n$th numerator is $zq^{{n\choose 2}}$. \vspace{.2in}
\end{theorem}

We think it is remarkable that such a  continued fraction encodes
information about $(132)$-avoiding permutations. We will first prove the
theorem, and then study some consequences and generalizations.
\section{The patterns}
Let the weight of a permutation $\pi$ of $|\pi|$ letters be
$z^{|\pi|}q^{|123(\pi)|}t^{|12(\pi)|}$, in which  $|123(\pi)|$ is the
number of patterns $(123)$ (rising triples) in $\pi$, and $|12(\pi)|$ is
the number of rising pairs in $\pi$. Let 
\begin{equation}
\label{eq:gfun}
P(q,z,t)= {\sum}'_{\pi}\mathrm{Weight}(\pi),
\end{equation}
where the sum extends over all $(132)$-avoiding permutations $\pi$. 

   If $\pi$ is a $(132)$-avoiding permutation on $\{1,2,\dots ,n\}$,
$(n>0)$ and 
 the largest element, $n$, is at the $k$th position, i.e., $\pi(k)=n$, then
by letting
 $\pi_1:=\{\pi(i)\}_1^{k-1}$ and $\pi_2:=\{\pi(i)\}_{k+1}^n$, we have that
every element in
 $\pi_1$ must be larger than every element of $\pi_2$,
 or else a $(132)$ would be formed, with the $n$ serving as the `$3$' of
the $(132)$.
 Hence, $\pi_1$ is a permutation of the set $\{n-k+1,\dots ,n-1\}$, and
$\pi_2$ is a permutation of the set  $\{1,\dots ,n-k\}$. Furthermore,
$\pi_1$ and $\pi_2$ are each $(132)$-avoiding. Conversely, if $\pi_1$ and
$\pi_2$ are $(132)$-avoiding permutations on
 $\{n-k+1,\dots ,n-1\}$ and $\{1,\dots ,n-k\}$ respectively (for some $k$,
 $1\le k\le n$), then $(\pi_1n\pi_2)$ is a nonempty $(132)$-avoiding
permutation.

 Thus we have:
\[        |123(\pi)| = |123(\pi_1)| + |123(\pi_2)| + |12(\pi_1)|,\] 
 since a $(123)$ pattern in $\pi\eqdef (\pi_1n\pi_2)$ may either be totally
immersed 
 in the $\pi_1$ part, or wholly immersed in the $\pi_2$ part, or may be due
to the $n$ 
 serving as the `$3$' of the $(123)$, the number of which is the number of
$(12)$  
 patterns in $\pi_1$.

 We also have 
\[        |12(\pi)| = |12(\pi_1)| + |12(\pi_2)| + |\pi_1|, \]
 and, of course
\[        |\pi| = |\pi_1| + |\pi_2| + 1.\]
 Hence,
\begin{eqnarray*}
 \mathrm{Weight}(\pi)(q,z,t)&:=&q^{|123(\pi)|}z^{|\pi|}t^{|12(\pi)|}\\
&=& q^{|123(\pi_1)|+|123(\pi_2)|+|12(\pi_1)|}z^{|\pi_1|+|\pi_2|+1}
 t^{|12(\pi_1)|+ |12(\pi_2)|+|\pi_1|}\\
&=&zq^{|123(\pi_1)|}(qt)^{|12(\pi_1)|}
 (zt)^{|\pi_1|} q^{|123(\pi_2)|}t^{|12(\pi_2)|}z^{|\pi_2|}\\
&=&
 z\mathrm{Weight}(\pi_1)(q,zt,tq)\mathrm{Weight}(\pi_2)(q,z,t).
\end{eqnarray*} 
 Now sum over all possible $(132)$-avoiding permutations $\pi$, to get the
functional equation 
\begin{equation}\label{eq:fcnleq}     P(q,z,t) = 1 +
zP(q,zt,tq)P(q,z,t),\end{equation}
in which the 1 corresponds to the empty permutation.

    Next let $Q(q,z,t)$ be the sum of all the weights of all permutations with 
 exactly ONE $(132)$ pattern.  By adapting the argument from Mikl\' os B\'
ona's paper \cite{Bo} 
 we easily see that $Q(q,z,t)$ satisfies
\begin{equation}\label{eq:qfcnleq}
Q(q,z,t) = zP(q,zt,qt)Q(q,z,t) + zQ(q,zt,qt)P(q,z,t) +
               t^2z^2P(q,zt,qt)(P(q,z,t)-1) .\end{equation}
    This holds since our sole $(132)$ pattern can either appear in the elements
\begin{enumerate} 
\item[(a)] before $n$, \item[(b)] after $n$, or \item[(c)] with $n$ as the
`3' in the $(132)$ pattern. 
\end{enumerate} 
 The term $zP(q,zt,qt)Q(q,z,t)$ corresponds to (a),  $zQ(q,zt,qt)P(q,z,t)$  
 corresponds to (b), and $t^2z^2P(q,zt,qt)(P(q,z,t)-1)$ corresponds to 
 (c). We see that case (c) follows since  $\pi=(\pi_1,n-k,n,\pi_2)$, where
 $\pi_1$ is a permutation of $[n-k+2,\dots ,n-1]$, $\pi_2$ is a permutation of 
 $[1,\dots ,n-k-1]\cup \{n-k+1\}$, and $k\neq n$.

\section{The fractions}
Here we study this generating function $P(q,z,t)$ further, and find that it
is a pretty continued fraction, and discover a fairly explicit form for its
numerator and denominator.

First, from (\ref{eq:fcnleq}) we have that
\begin{equation}\label{eq:cfrac}P(q,z,t)=\frac{1}{1-zP(q,zt,tq)},\end{equation}
and so by iteration we have the continued fraction,
\begin{equation}\label{eq:cfrac2}
P(q,z,t)={1\over\displaystyle 1-{\strut z\over\displaystyle 1-{\strut
zt\over\displaystyle 1-{\strut zt^2q\over\displaystyle 1-{\strut
zt^3q^3\over\displaystyle 1-{\strut zt^4q^6\over\displaystyle \dots}}}}}}
\end{equation}

Now let
\[P(q,z,t)=\frac{A(q,z,t)}{B(q,z,t)}.\]
Then substitution in (\ref{eq:cfrac}) shows that $A(q,z,t)=B(q,zt,tq)$, and
therefore
\begin{equation}
\label{eq:pform}
P(q,z,t)=\frac{B(q,zt,tq)}{B(q,z,t)}
\end{equation}
where $B$ satisfies the functional equation
\begin{equation}
\label{eq:beq}
B(q,z,t)=B(q,zt,tq)-zB(q,t^2qz,tq^2).
\end{equation}

To find out more about the form of $B$ we write
\[B(q,z,t)=\sum_{m\ge 0}\phi_m(q,t)z^m.\]
Then $\phi_0=1$, and
\[\phi_m(q,t)=t^m\phi_m(q,qt)-t^{2m-2}q^{m-1}\phi_{m-1}(q,tq^2),\]
for $m=1,2,3,\dots $. It is easy to see by induction that
\[\phi_m(q,t)=-\sum_{j\ge 2}t^{jm-2}q^{m{j\choose
2}-2j+3}\phi_{m-1}(q,tq^j).\qquad (m\ge 1;\phi_0=1)\]
For example, we have
\[\phi_1(q,t)=-\sum_{j\ge 0}t^jq^{{j\choose 2}},\]
and
\[\phi_2(q,t)=\sum_{j,\ell\ge 2}t^{\ell+2j-4}q^{\frac{1}{2}\ell^2+j^2+\ell
j-\frac{5}{2}\ell-5j+6}.\]

In general, the exponent of $t$ in $\phi_m(q,t)$ will be a linear form in
the $m$ summation indices, plus a constant, and the exponent of $q$ will be
an affine form in these indices, i.e., a quadratic form plus a linear form
plus a constant. Let's find all of these forms explicitly.

Hence suppose in general that
\[
\phi_m(t)=(-1)^m\sum_{\mathbf{j}\ge 0}t^{\mathbf{a}_m\cdot
\mathbf{j}+b_m}q^{(\mathbf{j},Q_m\mathbf{j})+\mathbf{c}_m\cdot\mathbf{j}+d_m},
\]
in which $\mathbf{j}$ is the $m$-vector of summation indices, $Q_m$ is a
real symmetric $m\times m$ matrix to be determined,
$\mathbf{a}_m,\mathbf{c}_m$ are $m$-vectors, and $b_m,d_m$ are scalars.
Inductively we find that
\begin{eqnarray*}\label{eq:avec}
\mathbf{a}_m&=&\{r\}_{r=1}^m,\\
b_m&=&-2m,\\
\mathbf{c}_m&=&\left\{-5r/2\right\}_{r=1}^m,\\
d_m&=&3m.\end{eqnarray*}
The $m\times m$ matrix $Q_m$ is $\{\min{(r,s)}/2\}_{r,s=1}^m$. Thus we have
the following formula for $B$.
\begin{theorem}
The denominator $B(q,z,t)$ of the grand generating function $P(q,z,t)$ is
explicitly given by
\begin{equation}
\label{eq:bform}
B(q,z,t)=1+\sum_{m=1}^{\infty}(-zq^3t^{-2})^m\sum_{j_1,\dots ,j_m\ge
2}t^{\sum_{r=1}^mrj_r}q^{\frac{1}{2}\left\{\sum_{r,s=1}^m\min{(r,s)}j_rj_s-5
\sum_{r=1}^mrj_r\right\}}.
\end{equation}
\end{theorem}

\section{The series computations}

If $f_r(n)$ denotes the number of permutations of $n$ letters that contain
no pattern $(132)$ and have exactly $r$ $(123)$'s, we write
AR$(r,z):=\sum_nf_r(n)z^n$. Then AR$(r,z)$ is the coefficient of $q^r$ in
the series development of $P(q,z,1)$ of (\ref{eq:gfun}). That is, we have
\begin{equation}
\label{eq:argen2}
{1\over\displaystyle 1-{\strut z\over\displaystyle 1-{\strut
z\over\displaystyle 1-{\strut zq\over\displaystyle 1-{\strut
zq^3\over\displaystyle 1-{\strut zq^6\over\displaystyle
\dots}}}}}}=\sum_{r\ge 0}\mathrm{AR}(r,z)q^r
\end{equation}
>From (\ref{eq:argen2}) we see that if we terminate the fraction $P(q,z,1)$
at the numerator $q^N$, say, then we'll know all of the
$\{\mathrm{AR}(r,z)\}_{r=0}^N$ exactly.

Further, if we know the denominator $B(q,z,t)$ in (\ref{eq:pform}) exactly
through terms of order $q^N$, then by carrying out the division in
(\ref{eq:pform}) and keeping the same accuracy, we will, after setting
$t=1$, again obtain all of the generating functions $\{$AR$(r,z)\}_{r=0}^N$
exactly.

Finally, to find the denominator $B(q,z,t)$ in (\ref{eq:pform}) exactly
through terms of order $q^N$, it is sufficient to carry out the iteration
that is implicit in (\ref{eq:beq}) $N$ times, since further iteration will
affect only the terms involving powers of $q$ higher than the $N$th.

In that way we computed the AR$(r,z)$'s for $0\le r\le 15$ in a few
seconds, as is shown below in the initial section of the series
(\ref{eq:argen2}):

\vspace{.2in}

\hspace{.1in}\mbox{${\frac{1-z}{1-2\,z}} + {\frac{{z^3}\,q}{{{\left( 1 -
2\,z \right) }^2}}} + 
  {\frac{\left( 1 - z \right) \,{z^4}\,{q^2}}{{{\left( 1 - 2\,z \right)
}^3}}} + 
  {\frac{{{\left( 1 - z \right) }^2}\,{z^5}\,{q^3}}{{{\left( 1 - 2\,z
\right) }^4}}} + 
  {\frac{{z^4}\,\left( -1 + 6\,z - 13\,{z^2} + 11\,{z^3} - 3\,{z^4} + {z^5}
\right) \,{q^4}}
    {{{\left( -1 + 2\,z \right) }^5}}} $}

\vspace{.1in}

\hspace{.3in}\mbox{$+ {\frac{{z^5}\,
      \left( 2 - 14\,z + 37\,{z^2} - 44\,{z^3} + 22\,{z^4} - 4\,{z^5} +
{z^6} \right) \,{q^5}}{{{\left( 1 - 
          2\,z \right) }^6}}}
 + {\frac{{{\left( 1 - z \right) }^2}\,{z^6}\,
      \left( -3 + 18\,z - 37\,{z^2} + 27\,{z^3} - 3\,{z^4} + {z^5} \right)
\,{q^6}}{{{\left( -1 + 2\,z \right) }^7}}} $}  

\vspace{.1in}

\hspace{.5in}\mbox{$+{\frac{{z^5}\,\left( 1 - 12\,z + 64\,{z^2} -
196\,{z^3} + 373\,{z^4} - 450\,{z^5} + 343\,{z^6} - 164\,{z^7} + 47\,{z^8}
- 6\,{z^9} + {z^{10}} \right) \,{q^7}}{{{\left( 1 - 2\,z \right)
}^8}}}+\dots $}

\vspace{.2in}

If $g_r(n)$ denotes the number of permutations of $n$ letters that
contain $1$ $(132)$ pattern and have exactly $r$ $(123)$'s, we 
write Aaron(r,z)$:=\sum_{n} g_r(n)z^n$.  Then Aaron(r,z) is 
the coefficient of $q^r$ in the series development of 
$Q(q,z,1)$ of (\ref{eq:qfcnleq}).  Since we have a very quick method to
compute $P(q,z,1)$, we can iterate equation (4) to
compute the Aaron(r,z)'s.  Shown below are the Aaron(r,z)'s
for $0 \leq r \leq 6$, that were computed in a few minutes.

\vspace{.2in}
$$
\begin{array}{l}
\frac{z^3}{(1-2z)^2} + \frac{2z^5}{(1-2z)^3}q
+ \frac{z^4(z^3-6z^2+4z-1)}{(1-2z)^4} q^2
+ \frac{2z^5(z-1)(5z^2-4z+1)}{(1-2z)^5} q^3\\
\\
+ \frac{z^6(z^5+12z^4-55z^3+65z^2-30z+5)}{(1-2z)^6} q^4
+ \frac{-2z^7(z^6+6z^5-40z^4+80z^3-69z^2+27z-4)}{(1-2z)^7} q^5\\
\\
+ \frac{-z^6(z-1)(3z^8+13z^7-77z^6+240z^5-329z^4+231z^3-91z^2+20z-2)}
          {(1-2z)^8} q^6
\end{array}
$$

\end{document}